\pdfoutput=1
 \documentclass[preprint,showpacs,preprintnumbers]{revtex4}
 
\usepackage{color} %,hyperref}
\usepackage{amsmath,amssymb,graphicx,amsthm}
\usepackage{times}
\usepackage[normalem]{ulem}
\usepackage{subfigure}

\def\ip<#1,#2>{\left\langle #1,#2\right\rangle}
\definecolor{brown}{rgb}{0.6,0.4,0.2}
  \def\todor#1{\relax\ignorespaces}
  \def\rev#1{{#1}}
  \def\todo#1{\relax\ignorespaces}
  \def\sout#1{\relax\ignorespaces}
  \def\remark#1{\relax\ignorespaces}
\begin{document}

%opening
\title{Modeling of Transitional Channel Flow Using Balanced Proper Orthogonal Decomposition}

\author{Milo\v{s} Ilak}
\email{milak@princeton.edu}
\affiliation{%
Department of Mechanical and Aerospace Engineering, \\
 Princeton University, Princeton, NJ 08544, USA
}%
\author{Clarence W. Rowley}
\email{cwrowley@princeton.edu}
\affiliation{%
Department of Mechanical and Aerospace Engineering, \\
 Princeton University, Princeton, NJ 08544, USA
}%

\pacs{47.27.nd,47.27.Cn,47.85.L-}
% 47.27.nd is channel flow
%47.27.Cn is transition to turbulence
%47.85.L- is flow control
\date{\today}

\newcommand{\image}{\operatorname{image}}
\newcommand{\Tr}{\operatorname{Tr}}
                                        
%uncomment the following line if using article, 
%for aps format it needs to be after the abstract
%\maketitle

\begin{abstract}
%\remark{Remarks about the changes are in red.}
%\ifMarkup
 % \rev{Actual changes in the text are in blue.}\\
%\fi

We study reduced-order models of three-dimensional perturbations in linearized channel flow using balanced proper orthogonal decomposition (BPOD). The models are obtained from three-dimensional simulations in physical space as opposed to the traditional single-wavenumber approach, and are therefore \rev{better} able to capture the effects of localized disturbances or localized actuators.  In order to assess the performance of the models, we consider the impulse response and frequency response, and variation of the Reynolds number as a model parameter. We show that the BPOD procedure yields models that capture the transient growth well at a low order, whereas standard POD does not capture the growth unless a considerably larger number of modes is included, and even then can be inaccurate. In the case of a localized actuator, we show that POD modes which are not energetically significant can be very important for capturing the energy growth. In addition, a comparison of the subspaces resulting from the two methods suggests that the use of a non-orthogonal projection with adjoint modes is most likely the main reason for the superior performance of BPOD.  We also demonstrate that for single-wavenumber perturbations, low-order BPOD models reproduce the dominant eigenvalues of the full system better than POD models of the same order. These features indicate that the simple, yet accurate BPOD models are a good candidate for developing model-based controllers for channel flow. 

\end{abstract}

%comment out the following if using article class
\maketitle

%\tableofcontents

\section{Introduction}

Many techniques for developing practical controllers for fluids require models of the system that are both tractable and that describe accurately the flow physics for the given flow regime.  One of the problems of great interest in flow control is drag reduction in shear flows. Drag increases drastically as flow transitions from laminar to turbulent, making turbulence suppression or inhibition of transition to turbulence promising strategies for drag reduction. Both open-loop and closed-loop strategies have been used in control of channel flow. Recently, Min et al~\cite{Min-06} have achieved sub-laminar drag through open-loop control using a traveling wave actuation on the channel walls. Lee et al~\cite{LeeCortelezzi-01} have reported drag reduction using closed-loop controllers developed for the linearized flow, while Joshi et al~\cite{JosSpeKim-97} and H\"ogberg et al~\cite{HogbergBewley-03} have demonstrated significant reduction in the energy of the perturbations to laminar flow using closed-loop controllers. Despite these recent successes, the adequate modeling of the fluid flow and the actuators in a framework useful for practical control design is still a challenge. In order to use model-based control strategies, one needs an accurate description of the system dynamics and the actuation, from a model of sufficiently low order to allow practical implementation.  The goal of this paper is to explore improved techniques for developing such reduced-order models, in the context of a transitional channel flow.

The POD/Galerkin method has been used extensively for reduced-order modeling of fluid problems~\cite{Sirovich-87,MoinMoser-89,Smith-thesis,HLB-96,Rowley-physd04}. In this method, an empirical basis of orthonormal eigenfunctions is obtained from experimental or simulation data, and the Navier-Stokes equations are projected onto this basis. For fluid problems, this basis is optimal in terms of capturing the energy of the flow: the most significant modes are the ones that carry most of the kinetic energy. Although this method is applicable to many different types of flows, and computationally tractable for very large data sets, it can result in inaccurate low-order models since the most energetic modes are not always the most dynamically significant ones, as for example in the case of acoustic modes in cavity oscillations~\cite{Rowley-physd04}. Improved POD models can in some cases be obtained by a careful selection of modes~\cite{Smith-thesis}, removal of symmetry from the data using Fourier or traveling modes~\cite{RowlMars-00,IlaRow-06} or inclusion of shift modes~\cite{Noack-03}, but these techniques are often ad hoc and typically require extensive fine tuning. In addition, there is no standard \sout{way of} \rev{method for} incorporating actuation into POD models, and if using these models for control design, one needs to be careful to respect the region of validity of the model in the presence of actuation~\cite{Tadmor2004tcst2vort}.

On the other hand, many tools from the control theory community have been developed for model reduction of linearized dynamics, often with a priori error bounds~\cite{ObinataAnderson}.  Here we focus on one of these methods, balanced truncation, introduced by Moore~\cite{Moore-81}, and described in detail in standard references (\cite{DullerudPag,ZhoDoy-98,ObinataAnderson}). Balanced truncation involves transforming a state-space system to a coordinate system where the states that respond most strongly to inputs (most controllable states) are also the states that have the most influence on future outputs (most observable states). In this \textit{balanced realization}, weakly observable and controllable states can be truncated to form reduced-order models that capture well the input-output behavior of the system. The required coordinate transformation is known as the balancing transformation.

The main difficulty with using balanced truncation for large systems is computational expense. In a typical fluids problem, the number of states is on the order of $10^5$ or more, and since computing balancing transformations involves solving Lyapunov equations for matrices of dimension $n\times n$  (where $n>10^5$ is the number of states) as well as correspondingly large eigenvalue problems, traditional approaches to balanced truncations may not be feasible.  Recently, snapshot-based methods for computing balanced truncation have been suggested, both to extend the concepts to nonlinear systems~\cite{LallMarsden-02} and to large linear systems~\cite{Rowley-ijbc05}. For systems with large numbers of outputs (e.g., if the output is the entire state), the snapshot-based method requires a large number of adjoint simulations, which still makes the problem computationally intractable. A procedure we call Balanced POD (BPOD)~\cite{Rowley-ijbc05}, uses an output projection to reduce the number of necessary adjoint simulations, and has been shown to be a computationally feasible approximation to balanced truncation for examples of two- and three-dimensional perturbations to laminar channel flow~\cite{Rowley-ijbc05,IlaRow-06}.

Many recent works on transition in shear flows have focused on the large non-normal transient growth of exponentially stable linear perturbations to the laminar flow, which is thought to lead to the so-called `subcritical' or `bypass' transition~\cite{SchmidHenn,Trefethen-93,Farrell-88,ButlerFarrell-92,Bamieh-01,JovBam-05,Reddy-98}. A comprehensive treatment of the subject is given by Schmid and Henningson~\cite{SchmidHenn}. Suppression of this large growth is of interest in \sout{controls} \rev{control} applications for drag reduction.
Control and estimation of linearized channel flow was studied by H\"ogberg et al~\cite{HogbergBewley-03}. More recently, {\AA}kervik et al~\cite{Aakervik_et_al-jfm-07} \remark{reference updated from conference proceedings to a journal article} studied control of a cavity flow using global eigenmodes. In the works by Farrell and Ioannou~\cite{Farrell-01} and Lee et al~\cite{LeeCortelezzi-01}, balanced truncation was applied to linearized channel flow at particular wavenumbers where the standard algorithms are applicable, since the full system is one-dimensional.

The contribution of this work is the application of BPOD to transitional channel flow for localized disturbances, without modeling on a wavenumber-by-wavenumber basis. Although a lot has been learned from studying `canonical' classes of perturbations (streamwise vortices, oblique waves, Tollmien-Schlichting waves, simple localized disturbances~\cite{HennLundJoh-93}), the standard wavenumber-by-wavenumber analysis has its limitations when efficient implementation of closed-loop control is desired, \rev{since model reduction and control design would need to be performed at each wavenumber}. The perturbations that arise in real flows often have complex three-dimensional structure with contributions at a wide range of wavenumbers as well as non-periodic or stochastic components. More importantly, for modeling of realistic actuation devices the wavenumber approach would \sout{not be feasible} \rev{be very complex} except for the special case of actuators that excite only specific wavenumbers. It is therefore of interest to have a method where disturbances with complex structure are modeled without Fourier decomposition, in particular for controls applications where physically realizable localized actuators need to be considered. To the best of our knowledge, such low-order balanced truncation models have not been reported in literature. \rev{Although some success has been achieved for channel flow using control at all wavenumbers, the method we propose is able to extract dominant dynamics of the flow in physical space, resulting in simpler models. Moreover, this method for balanced truncation would be very advantageous for more complex geometries, such as spatially developing boundary layers.} 

We first obtain models of a single-wavenumber perturbation computed by Farrell~\cite{ButlerFarrell-92} in order to validate our numerical methods, and then apply BPOD to a localized body-force actuator. The most desirable features of a reduced-order model are close approximation of the dynamics of the original system, inclusion of actuation and validity over a wide range of parameters, and we show that our BPOD models have those features. We also show that BPOD models capture the effects of actuation better than standard POD and that the considerable improvement in the capturing of the dynamics for BPOD is due to the non-orthogonal projection used. 

The rest of this paper is organized as follows. In Section~\ref{sec:modelred}, we briefly overview the two model reduction methods we use. In Section~\ref{sec:chanflow}, we describe the governing equations and our choice of inner product for the adjoint simulations of the system. In Sections~\ref{sec:results-single} and~\ref{sec:results-3D}, we present low-order models for a single-wavenumber optimal perturbation and a localized actuator and discuss their performance. \rev{Finally}, in Section~\ref{sec:conclusions}, we describe our conclusions and directions for future work.

\section{Model Reduction via POD and Balanced POD}\label{sec:modelred}

\subsection{Proper Orthogonal Decomposition (POD)}
\label{sec:modelred-POD}

We first give a brief overview of the POD/Galerkin method; details can be found in standard references~\cite{Sirovich-87,HLB-96}. The idea of Galerkin projection is, given a system
\begin{equation}
\dot x=f(x),\qquad x(t)\in \mathcal{X},
\end{equation}
where $\mathcal{X}$ is a high-dimensional Hilbert space, to project onto a low-dimensional subspace $S\subset\mathcal{X}$. Proper Orthogonal Decomposition determines an orthogonal basis for such a subspace, which is obtained by solving the eigenvalue problem for $XX^T$ where $X$ is a matrix whose columns are simulation snapshots $x(t_k)$ at some times $t_k$. In this basis, we can represent the dynamics of $x(t)$ as
\begin{equation}
x(t) = \sum_{j=1}^m a_j(t) \theta_j,
\label{eq:PODsum}
\end{equation}
\remark{the $r(t)$ in the previous equation was replaced by $x(t)$ since this is still the full dynamics, just in a different basis} where $\theta_j$ are the time-independent basis functions (POD modes) and $a_j(t)$ are the corresponding time coefficients, which are obtained from
\begin{equation}
\dot{a}_j = \ip<\theta_j,f(r)>.
\label{eq:PODtimecoeffs}
\end{equation}
A reduced-order model of order~$r$ can be obtained as a set of ODEs for the time evolution of these coeffecients by projecting the original system onto the most significant POD modes in the system (i.e. including only the first $r$ modes where $r<m$). The method is applicable to both linear and nonlinear systems. For linear systems, the POD modes of data arising from the input-state impulse response are the most controllable modes of the linear system~\cite{Rowley-ijbc05}. However, both controllability and observability are important for the input-output behavior of a system, and POD often fails to capture highly observable modes. On the other hand, balanced truncation does take into account both of these properties, and we describe this method next. 

\subsection{Balanced POD}
\label{sec:modelred-OP-bpod}

Balanced truncation is a standard model reduction method~\cite{Moore-81,DullerudPag,ZhoDoy-98} used for stable linear input-output systems \sout{given by} \rev{of the form}
 \begin{equation}
 \begin{aligned}
 \dot x &= Ax + Bu\\
 y &= Cx
 \end{aligned}
 \label{eq:primal}
 \end{equation}
where $u\in\mathcal{U} = \mathbb{R}^p$ is the vector of inputs, $y\in\mathcal{Y} = \mathbb{R}^q$ is the output, $x\in\mathcal{X}=\mathbb{R}^n$ is the state vector \rev{(although in general all three spaces can be complex as well)}, and $A$, $B$, and $C$ are matrices of appropriate dimension. The idea of balancing is to find a change of coordinates in which the controllability and observability Gramians, defined by  %
 \begin{equation}
 W_c = \int_0^\infty e^{At}BB^+e^{A^+t}\,dt,\quad 
 W_o = \int_0^\infty e^{A^+t}C^+Ce^{At}\,dt,
 \label{eq:gramdef} 
 \end{equation}
are equal and diagonal. Here $A^+$, $B^+$ and $C^+$ define the corresponding adjoint system. It should be noted that in general $A^+ \neq A^T$, the two being equal only when the inner product used to derive the adjoint does not have an associated weight. It can be shown that balanced truncation does not depend on the choice of the inner product on the state space $\mathcal{X}$ (see Appendix~\ref{sec:App-indep_ip}), although it does depend on the choices of inner products for $\mathcal{U}$ and $\mathcal{Y}$ . One then truncates the least controllable and observable modes, corresponding to the smallest eigenvalues of these Gramians.  A detailed description of the Balanced POD procedure, which is a computationally tractable procedure for finding such a transformation, is given in~\cite{Rowley-ijbc05}. In this method, one begins by computing snapshots of the impulse-state response of the system in Eq.~(\ref{eq:primal}) and the adjoint system
\begin{equation}
	\dot z = A^+z + C^+ v
	\label{eq:adjoint_sys}
\end{equation}
and stacking the direct and adjoint snapshots as columns of matrices $X$ and~$Y$ (with appropriate quadrature weights~\cite{Rowley-ijbc05}).  One can show that the Gramians in Eq.~(\ref{eq:gramdef}) may then be approximated by \textit{empirical Gramians}~\cite{LallMarsden-02} $W_{c,e}$ and~$W_{o,e}$, as
\begin{equation}
	W_c \approx W_{c,e} = XX^+,\qquad W_o \approx W_{o,e}=YY^+.
	\label{eq:empgrams}
\end{equation}
The key idea in the method of snapshots is to compute the transformation that balances the empirical Gramians (or at least the dominant directions of this transformation) without actually computing the Gramians themselves, whose dimension is large.  In this respect, the method of snapshots for BPOD~\cite{Rowley-ijbc05} resembles the method of snapshots introduced by Sirovich~\cite{Sirovich-87} for more efficient computation of POD modes.  To compute the balancing transformation, one computes the singular value decomposition (SVD) of the matrix $Y^+X$ \rev{(see Appendix~\ref{sec:App-indep_ip} for a discussion of $Y^+$)}:
 \begin{equation}
Y^+X = U\Sigma V^T,
\label{eq:SVD}
\end{equation}
from which the balancing transformation $\Phi$ and its inverse $\Psi$ are found by
\begin{equation}
\Phi = XV\Sigma^{-1/2},\quad \Psi = YU\Sigma^{-1/2}.
\label{eq:ST_def}
\end{equation}
The columns of $\Phi$ are the \textit{balancing modes} and the columns of $\Psi$ are the \textit{adjoint modes}, and the two sets of modes are biorthogonal. The entries of the diagonal matrix $\Sigma$ are known as the Hankel singular values (HSVs). The non-orthogonal projection onto the basis of balancing modes using the adjoint modes is also known as Petrov-Galerkin projection.

Note that a different procedure for approximating balancing transformations has also been used in~\cite{Willcox-02}, in which the Gramians are separately reduced (that is, low-rank approximations of $W_{c,e}$ and $W_{o,e}$ are first constructed, and then the balancing transformation for the rank-reduced Gramians is computed by an unspecified algorithm).  However, this procedure is more computationally intensive than our procedure, and also gives worse results~\cite{Rowley-ijbc05}, since almost-uncontrollable modes may be strongly observable, so should not be truncated.

\subsection{Output projection}
If the number of outputs of the system is large, as in a typical fluids problem (where $n=q$, i.e. the output is the full state), the computation of the adjoint simulations of the system given by Eq.~(\ref{eq:adjoint_sys}) may not be tractable, since one simulation is needed for each component of the output. 
A way to reduce the number of system outputs is to first project the output onto a low-dimensional subspace, i.e., taking $\tilde{y}=P_{s} Cx$, where $P_{s}$ is an orthogonal projection onto a $s$-dimensional subspace of~$\mathcal{Y}$, as suggested in~\cite{Rowley-ijbc05}.  The system is now of the form:
\begin{equation}
\begin{aligned}
\dot x &= Ax + Bu\\
\tilde{y} &= P_{s}Cx
\end{aligned}
\label{eq:OPprimal}
\end{equation}
where $s$ is the rank of the output projection. The projection~$P_{s}$ that minimizes the 2-norm of the difference between the original transfer function and the output-projected transfer function is given simply by the POD of the set of impulse-state responses~\cite{Rowley-ijbc05}. This projection can be written as $P_{s}=\Theta_{s}\Theta_{s}^+$, 
where columns of $\Theta_{s}:\mathbb{R}^{s}\to\mathcal{Y}$ are POD modes. Another way to write the system is as follows:
\begin{equation}
\begin{aligned}
\dot x &= Ax + Bu\\
\hat{y} &= \Theta_{s}^+Cx
\end{aligned}
\label{eq:OPprimal_coefs}
\end{equation}
Here, the outputs of the system are just the coefficients of the POD modes of the system impulse response and $\hat{y} \in \mathbb{R}^s$. This $s$-dimensional output carries the same information as the $n$-dimensional output $\tilde{y}$, which can be shown by Parseval's theorem. The corresponding adjoint system can now be written as
\begin{equation}
\begin{aligned}
\dot z &= A^+z + (\Theta_{s}^+C)^+v\\
w &= B^+z.
\end{aligned}
\label{eq:OPadj2}
\end{equation}
Note that if the output is the full state, so that $C=I$, and the adjoint is defined with respect to the standard $L_2$ inner product, the initial conditions of the adjoint simulations are just the POD modes (columns of $\Theta_{s}$). In practical computations, depending on the choice of inner product used in defining the adjoint system, and on the numerical quadrature method (for example, if the computations are done using Chebyshev polynomials) the matrix $\Theta_{s}^+$ is usually just the matrix $\Theta_{s}$ pre-multiplied by a matrix of inner product weights.

The idea of Balanced POD is to compute the snapshot-based balanced truncation of the system in Eq.~(\ref{eq:OPprimal_coefs}) instead of Eq.~(\ref{eq:primal}), so that only $s$ adjoint simulations are needed. It is easily shown that the systems in Eqs.~(\ref{eq:OPprimal}) and~(\ref{eq:OPprimal_coefs}) have the same observability Gramian using the fact that for any projection $P$, we have $P^2=P$. Transforming Eq.~(\ref{eq:OPprimal_coefs}) to balanced coordinates and writing $x=\Phi_1a$ is obtained as follows:
\begin{equation}\label{eq:balanced_system}
\begin{aligned}
\dot{a} &=\Psi_1^+A\Phi_1a+\Psi_1^+ Bu\\
y &=\Theta_{s}^+C\Phi_1a.
\end{aligned}
\end{equation}
The inverse transformation matrix $\Psi_{1}:\mathbb{R}^{r}\to\mathcal{X}$ and the transformation matrix $\Phi_{1}:\mathbb{R}^{r}\to\mathcal{X}$ are $n \times r$, $r$ being the number of states we want to retain in the system, which we will refer to as the rank of the model. Note that $r\le p$, where $p$ is the number of non-zero HSVs. For simplicity, we will assume from now on that $C=I_n$, i.e. the output of the original system is the full state \rev{(this is the case in fluid simulations in which we need to know the entire flow field)}. We can then represent the output of Eq.~(\ref{eq:balanced_system}) as $y=\Theta_{s}^+\Phi_1 z$, which is now the vector of time coefficients of the ${s}$ standard POD modes obtained from the impulse response of the system. For fluids systems the full flow field output of the model can be recovered from these coefficients and the corresponding modes. For a given dimension of the output projection, all BPOD models will have $s$ outputs regardless of the model rank $r$, while the number of POD model outputs is equal to $r$ at each rank. The effect of output projection on model performance will be illustrated in Sec.~\ref{sec:results-single} and Sec.~\ref{sec:results-3D}.
\section{Application to Transitional Channel Flow}\label{sec:chanflow}

\subsection{Governing equations}
\label{sec:chanflow-eqns}

For shear flows, the linearized equations may be conveniently written in terms of the wall-normal velocity~$v$ and the wall-normal vorticity~$\eta$ (see, for instance, \cite{SchmidHenn}).  The other variables (e.g., streamwise and spanwise velocities $u$ and $w$) may then be computed using the continuity equation $\partial_x u + \partial_y v + \partial_z w=0$ and the definition of wall-normal vorticity.  In these coordinates, the linearized (nondimensional) equations have the form
\begin{eqnarray}\label{eq:OSandSQeqns}
\left[(\partial_t + U\partial_x)\Delta-U''\partial_x -\frac{1}{Re}\Delta^2\right]v &=&0 \\
\left[\partial_t+U\partial_ x-\frac{1}{Re}\Delta\right]\eta&=&-U'\frac{\partial v}{\partial z}.
\end{eqnarray}
Here, $Re=U_{c}\delta/\nu$ is the Reynolds number, where $\nu$ is the kinematic viscosity, $\delta$ is the half-width of the channel, and $\Delta=\partial_x^2 + \partial_y^2 + \partial_z^2$ is the Laplacian. $U_c$ is a characteristic velocity, which for linearized channel flow is the centerline velocity of the laminar profile $U(y)$. The prime indicates differentiation with respect to $y$. The first equation is the Orr-Sommerfeld equation and the second one is known as the Squire equation. It was first shown numerically by Orszag~\cite{Orszag-1971} that the Orr-Sommerfeld equation for channel flow is stable up to $Re\approx 5772$, when an exponentially unstable eigenmode first arises. The Squire equation has stable eigenmodes for all values of $Re$. Still, complex behavior due to the non-normality exists for stable eigenmodes. The term on the right hand side of the Squire equation represents tilting of the spanwise component of the vorticity of the mean flow (which here is just $U'$) by the strain rate $\partial v/\partial z$~\cite{ButlerFarrell-92}, which gives rise to wall-normal vorticity. In the limit of high Reynolds number, the perturbation growth is dominated by this process, in particular for streamwise-constant perturbations. While the system also exhibits phenomena such as degeneracies and resonances~\cite{Gustavsson-86,HennSchmid-92}, non-normality has been shown to have a dominating effect on the energy growth~\cite{ReddyHenn-93}.

In operator form, we can represent the equations using more compact notation as follows:

\begin{equation}
\frac{\partial}{\partial t}
\begin{bmatrix}
-\Delta&0\\0&I
\end{bmatrix}
\begin{bmatrix} v\\ \eta\end{bmatrix}
= 
\begin{bmatrix}
L_{OS} & 0\\ -U'\partial_z & L_{SQ}
\end{bmatrix}
\begin{bmatrix} v\\ \eta\end{bmatrix}
\label{eq:linearized}
\end{equation}
where 
\begin{align*}
L_{OS} &= U\partial_x\Delta -U''\partial_x - \frac{1}{Re}\Delta^2\\
L_{SQ} &= -U\partial_x + \frac{1}{Re}\Delta
\end{align*}
are the Orr-Sommerfeld and Squire operators, respectively. If we define 
\begin{equation}
A=
\begin{bmatrix}
-\Delta&0\\0&I
\end{bmatrix}^{-1}
\begin{bmatrix}
L_{OS} & 0\\ -U'\partial_z & L_{SQ}
\end{bmatrix}
\label{eq:Amat_def}
\end{equation}
with no-slip boundary conditions, we can write the system in standard state-space form:
\begin{equation}
\begin{aligned}
\dot x &= Ax + Bu_1 + Fu_2\\
y &= Cx
\end{aligned}
\label{eq:ss_act_dist}
\end{equation}
where \sout{$B$ represents actuation and $F$ represents disturbances} \rev{$B$ and $F$ represent the spatial distributions of the actuators and disturbances respectively}, with $u_1(t)$ and $u_2(t)$ being the corresponding input vectors \rev{(the time-dependent amplitudes of the columns of $B$ and $F$)}. \rev{The actuation and the disturbances are equivalent mathematically as they are both inputs to the system.} We note here that the impulse-state responses are given by $x_1(t) = e^{At}B$ and $x_2(t) = e^{At}F$, and the adjoint system impulse-state responses for the full system are given by $z(t)=e^{A^+t}C^+$. Therefore, to obtain the POD basis needed for BPOD, we simulate the system given by Eq.~(\ref{eq:OSandSQeqns}) with a given perturbation or actuator as initial condition until the response has decayed to negligible levels, so that the matrix $XX^T$ that can be formed from the snapshots will closely approximate the controllability Gramian given by Eq.~(\ref{eq:gramdef}), where the integral extends to infinite time. \rev{Of course, computation of the matrix $XX^T$ is intractable for very large systems, so we compute} \sout{For} POD \sout{computation} via the method of snapshots, \sout{we form} \rev{forming} the smaller matrix $X^TX$ and follow\rev{ing} the procedure described in Sec.~\ref{sec:modelred-POD}. 

\subsection{Inner product on the state space}
\label{sec:chanflow-ip}

To determine the corresponding adjoint equations, one first needs to define an inner product on the vector space $\mathcal{X}$ of flow variables $(v,\eta)$. Since balanced truncation is independent of the choice of inner product used to define the adjoint (see Appendix~\ref{sec:App-indep_ip}), we may choose an inner product which is convenient for numerical computations. Let us define the inner product
\begin{equation}
\ip<(v_1,\eta_1),(v_2,\eta_2)>_M = \int_\Omega (-v_1\Delta v_2 + \eta_1\eta_2)\,dx\,dy\,dz,
\label{eq:define_Mip}
\end{equation}
where $\Omega$ denotes the fluid volume.  Note that, letting $M:\mathcal{X}\to\mathcal{X}$ denote the matrix operator on the left hand side of Eq.~(\ref{eq:linearized}), this is just the $L_2$ inner product of $(v_1,\eta_1)$ with $M(v_2,\eta_2)$. This inner product is different from the standard energy inner product used in analyzing perturbations through Fourier decomposition~\cite{Gustavsson-86,ButlerFarrell-92}, as there is no re-scaling at each wavenumber. 

With this definition of the inner product, the adjoint equations are easily found by integration by parts:
\begin{equation}
\frac{\partial}{\partial t}
\begin{bmatrix}
-\Delta&0\\0&I
\end{bmatrix}
\begin{bmatrix} v\\ \eta\end{bmatrix}
= 
\begin{bmatrix}
L_{OS}^* & U'\partial_z\\ 0 & L_{SQ}^*
\end{bmatrix}
\begin{bmatrix} v\\ \eta\end{bmatrix}
\label{eq:adjoint}
\end{equation}
where
\begin{align*}
L_{OS}^* &= -U\partial_x\Delta -2 U'\partial_x\partial_y - \frac{1}{Re}\Delta^2\\
L_{SQ}^* &= U\partial_x + \frac{1}{Re}\Delta.
\end{align*}

\subsection{Inner product on the output space}\label{sec:results-ipchoice}

Although the time evolution of the linearized disturbances is fully described by the wall-normal velocity-vorticity formulation, the output of the system can be chosen to be in different variables. When using POD, the choice of inner product can have a large impact on the results. If the output of our system is only the velocity-vorticity field, the standard $L_2$ inner product can be used. For our system, since the other two velocity components can easily be recovered using continuity and the definition of vorticity, we can choose the full velocity field to be the output, and use the energy inner product given by
\begin{equation}
\left\langle\mathbf{u_1},\mathbf{u_2}\right\rangle = \int_\Omega (u_1u_2 + v_1v_2 + w_1w_2)\,dx\,dy\,dz,
\label{eq:ip_energy}
\end{equation}
%
%\remark{Milos, I removed parentheses from the LHS of the above equation.}  
%
This choice is more intuitively appealing, since the POD modes for the output projection will capture the true kinetic energy of the perturbation. We therefore define the output space $\mathcal{Y}$ in our system as the space $\mathbb{R}^n$ together with the inner product defined by Eq.~(\ref{eq:ip_energy}). We note here that the space $\mathcal{X}$ is also $\mathbb{R}^n$, though endowed with a different inner product (the $M$-inner product described in the previous section).

\subsection{Numerical Methods}

The simulations were performed using a linearized version of a fully nonlinear DNS code using the spectral \sout{collocation} method described by Kim et~al.~\cite{KimMoinMoser-87}, with periodic boundary conditions in the streamwise and spanwise directions. The linearized code was verified against the analytic time evolution of Orr-Sommerfeld eigenfunctions and optimal perturbations, and the resolution of each simulation was checked by varying the time step and grid resolution. The size of the computational box was $2\pi \times 2\pi$ in the streamwise and spanwise directions for all simulations. Standard LAPACK routines were used for the computation of POD and balanced POD modes, as well as for the comparison of subspaces. The reduced-order models were integrated using the standard fourth-order Runge-Kutta scheme. All computations were done using Fortran 90 and MATLAB. A code by S.C.Reddy~\cite{SchmidHenn} was used to compute the initial conditions for the single-wavenumber perturbations. The integration weights derived by Hanifi et al~\cite{Hanifi-96} were used for the computation of inner products on the Chebyshev grid.

\section{Results: single-wavenumber perturbations}\label{sec:results-single}

\remark{The first few paragraphs of this section have some considerable additions to the text to address the concern of the first reviewer about the single wavenumber results.}% Hopefully it is now clearly explained.}

We start by investigating the system given by Eq.~(\ref{eq:ss_act_dist}) without actuation and \rev{only} in the presence of disturbances \rev{(without the $Bu_1$ term)}. In order to validate the numerical methods, we first obtain BPOD models from three-dimensional simulations of simple and well-known single-wavenumber perturbation cases, described by Butler and Farrell~\cite{ButlerFarrell-92} and also investigated by Schmid and Henningson~\cite{SchmidHenn}. The general form of such disturbances is given by
\begin{equation}\label{eq:1Dpert_def}
q(x,y,z,t) = \hat{q}(y,t)e^{(i\alpha x+i\beta z)},
\end{equation}
where $\hat{q}(y,t)=[v(y,t)\,\,\eta(y,t)]^T$. The standard approach to such perturbations is to compute the time evolution of $\hat{q}(y,t)$, which fully describes the system, since the velocity components $u$ and $w$ can easily be computed. For this one-dimensional problem, standard algorithms for computing balanced truncation are computationally tractable. Therefore, we are able to compare the models resulting from exact balanced truncation \rev{(which for the 1-D case can be computed in MATLAB using standard algorithms~\cite{Laub-87})} to BPOD models obtained from three-dimensional simulations of the real part of the full field, $Re\left\lbrace q(x,y,z,t)\right\rbrace $ \rev{at a particular wavenumber pair ($\alpha,\beta$)} on a large grid, similar to the comparison made by Rowley~\cite{Rowley-ijbc05} for a streamwise-constant perturbation. \rev{We note that for a given wavenumber pair the comparison between BPOD and exact balanced truncation can be done only using 1-D simulations, but we also performed 3-D simulations in order to verify our codes.}
% It can be shown that the results from the 3-D simulation are equal to the results of a 1-D system multiplied by a factor of $2\pi^2$. } \todor{Milos, what does this mean that the ``results'' are multiplied by $2\pi^2$?  What exactly is multiplied?  The velocities?  The energies?  HSVs?  Is $2\pi^2$ the volume, or why this factor?}  
\rev{We also note that, since the outputs of the output-projected system and the reduced-order models are coefficients of POD modes, the $C$ matrix in~(\ref{eq:ss_act_dist}) was modified so that the output of the full system is in the POD basis as well. This way, a meaningful comparison between the balanced truncation of the full system and BPOD is obtained.}

\begin{figure*}%[htpb]
 \centering
 \subfigure[]{\includegraphics[width=0.48\linewidth]{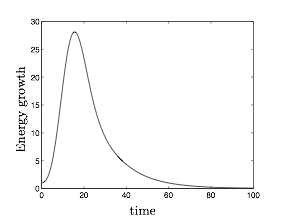}}
 \subfigure[]{\includegraphics[width=0.48\linewidth]{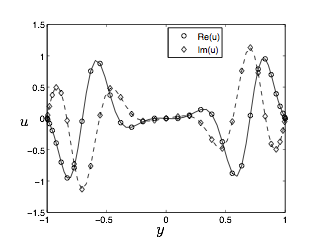}}
 \caption{(a) Kinetic energy growth for the optimal perturbation at wavenumber $\alpha=1,\beta=1$ at $Re=1000$. (b) The $\alpha=1,\beta=1$ optimal perturbation, showing streamwise velocity $u$ (complex).}
 \label{fig:en_growth_singlewavenum}
 \end{figure*}

The initial conditions were computed using the method described by Reddy and Henningson~\cite{ReddyHenn-93} and their energy growth was verified against values reported in that work. While streamwise-constant perturbations exhibit the largest energy growth, three-dimensional perturbations exhibit more interesting dynamics. \rev{(Here, by ``three-dimensional,'' we mean that the perturbations have components in both streamwise and spanwise directions, although the problem can still be treated as 1-D in the wall-normal direction as described above.)} We focus on the $\alpha=1,\beta=1$ perturbation at $Re=1000$, whose energy growth is shown in Fig.~\ref{fig:en_growth_singlewavenum}. The computational grid used in the three-dimensional simulation was $16 \times 65 \times 16$, corresponding to 33280 states in the system given by Eq.~(\ref{eq:ss_act_dist}). Balanced truncation of the one-dimensional problem with $65$ Chebyshev modes is easily and accurately computed using the algorithm described in~\cite{Rowley-ijbc05} so that BPOD performed on the large system can be compared to exact balanced truncation.

\subsection{Mode subspaces}\label{sec:mode_subspaces_single}

 \begin{figure*}%[htpb]
 \centering
 \subfigure[]{\includegraphics[width=0.48\linewidth]{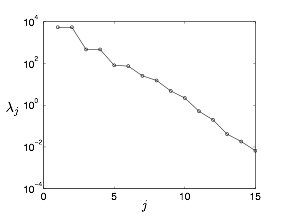}}
 \subfigure[]{\includegraphics[width=0.48\linewidth]{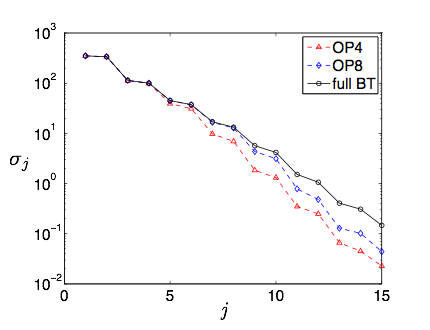}}
 \caption{(a) The first 15 POD eigenvalues for $\alpha=1,\beta=1$ initial perturbation  at $Re=1000$.
 (b) The first 15 Hankel singular values (HSVs) for: four-mode ($\triangle$) and eight-mode ($\diamond$)  output projections and full balanced truncation ($\circ$) for the same case.}
 \label{fig:a1b1_spectra}
 \end{figure*}

It was found that 500 equally spaced snapshots are sufficient for accurate computation of the POD modes, since for a larger number of snapshots with finer spacing there is no considerable change in the eigenvalue spectrum or the corresponding modes. 
We see from Fig.~\ref{fig:a1b1_spectra} that the most significant eigenvalues and the corresponding modes typically come in pairs, representing traveling structures that are 90 degrees out of phase. The first pair of modes contains $90.45\%$ of the energy, while the first three pairs contain $99.6\%$ of the energy. For the balanced POD models, a four-mode and eight-mode output projections were chosen, corresponding to respectively $98.3\%$ and $99.9\%$ of total energy contained in the POD modes.

\begin{figure}%[htpb]
\centering
\includegraphics[width=0.50\linewidth]{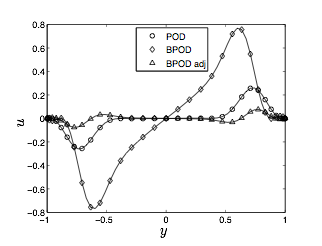}
\caption{Streamwise velocity for the first POD mode, balancing mode and adjoint mode for the $\alpha=1,\beta=1$, $Re=1000$ initial condition.}
\label{fig:a1b1_modes}
\end{figure}

We also notice that the HSVs (Fig.~\ref{fig:a1b1_spectra}) come in pairs, indicating that the most significant modes in the BPOD mode basis are again traveling structures similar to the standard POD modes. It is important to include these pairs of modes in the reduced-order models, as stability of the models for balanced truncation is guaranteed only if $\sigma_r>\sigma_{r+1}$ where $r$ is the rank of the model~\cite{DullerudPag}. While for standard POD modes there is no such requirement, mode pairs should always be included in the models on physical grounds. We also notice that the number of HSVs for each output projection that is equal to the full balanced truncation HSVs is approximately equal to the output projection rank. The same observation was made by Rowley~\cite{Rowley-ijbc05}, although there is no proof of this property at this point.

\begin{figure*}%[htpb]
\centering
\subfigure[]{\includegraphics[width=0.48\linewidth]{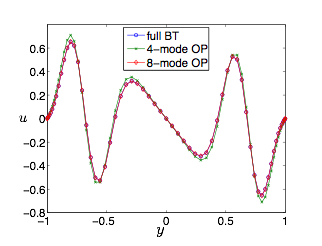}}
\subfigure[]{\includegraphics[width=0.48\linewidth]{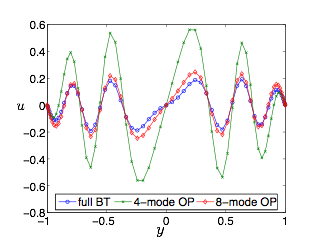}}
\caption{Streamwise velocity for (a) the sixth balancing mode and (b) the tenth balancing mode for BPOD with two different output projections and for full balanced truncation.}
\label{fig:a1b1_balmodes}
\end{figure*}

The first POD mode is shown in Fig.~\ref{fig:a1b1_modes} together with the first balancing and adjoint modes from a four-mode output projection. Figure \ref{fig:a1b1_balmodes} shows the streamwise velocity of the sixth and tenth balancing modes, illustrating the effect of the choice of output projection rank. The first four balancing modes from BPOD look identical for both output projections, while the sixth mode is not very accurately captured by a four-mode output projection. Both output projections do not capture very accurately the higher modes such as the tenth mode, which is also illustrated by the HSVs in Fig.~\ref{fig:a1b1_spectra}. As we show below, this inaccuracy does not significantly affect model performance, since these higher modes are not very significant dynamically.

\begin{figure*}
\centering
\subfigure[]{\includegraphics[width=0.30\linewidth]{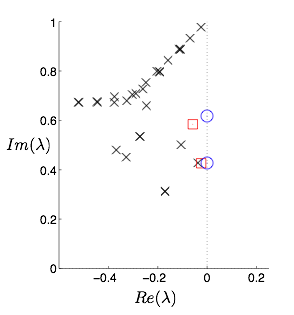}}
\subfigure[]{\includegraphics[width=0.30\linewidth]{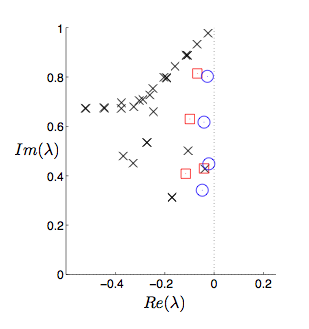}}
\subfigure[]{\includegraphics[width=0.30\linewidth]{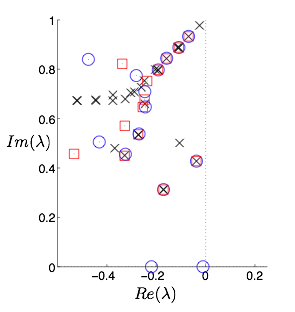}}
\caption{Spectrum of the full operator and reduced-order models for rank (a) 4, (b) 8,  (c) 30. The BPOD modes are from the eight-mode output projection. Symbols: \textcolor{red}{BPOD ($\square$)}, \textcolor{blue}{POD ($\bigcirc$)}, full operator ($+$) Only the most important part of the full spectrum is shown.}
\label{fig:a1b1_spec}
\end{figure*}

In this single-wavenumber case, the exact eigenvalue spectrum of the $A$ matrix from Eq.~(\ref{eq:Amat_def}) at a given Reynolds number can easily be computed. We note here that the eigenvalues of the matrix $A$ and therefore the poles of the corresponding transfer function are independent of the initial condition (which is just the $B$ matrix for our impulse response simulations). Figure~\ref{fig:a1b1_spec} shows the \sout{spectrum} \rev{spectra} of the full operator and three reduced-order models of different rank for the $\alpha=1,\beta=1$ perturbation. Since the \sout{spectrum is} \rev{spectra are} symmetric about the real axis, we only show the upper half of the complex plane. We see that, while the representation of the full spectrum improves for both methods as the rank increases, BPOD captures more accurately some of the most slowly decaying eigenvalues, which have the most influence on the dynamics of the system. For the rank four model, the standard POD model appears to be marginally stable, while the BPOD model approximates closely the eigenvalue closest to the origin. At higher order, the standard POD models improve and capture approximately the same eigenvalues as the BPOD models of the same rank. It is also important to notice that some of the eigenvalues of the full system are never captured by the reduced-order models. \rev{These eigenvalues correspond to uncontrollable eigenmodes of the full system, and can never be excited by this particular perturbation}.

\subsection{Impulse response}
\label{sec:results-singleimp}

\begin{figure*}%[htpb]
  \centering
\subfigure[]{\includegraphics[width=0.45\linewidth]{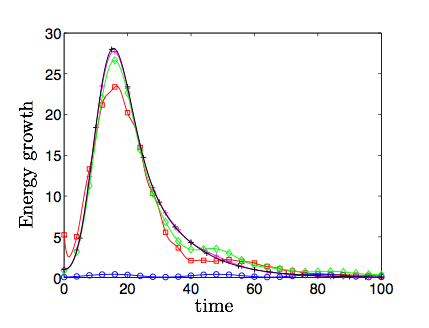}}
\subfigure[]{\includegraphics[width=0.45\linewidth]{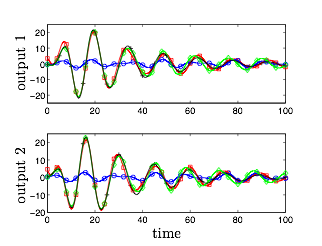}}
\caption{(a) $\alpha=1,\beta=1$ optimal perturbation at $Re=1000$, eight-mode output projection, 4-mode and 8-mode models. Full simulation ($+$), \textcolor{blue}{4-mode POD ($\bigcirc$)}, \textcolor{red}{4-mode BPOD ($\square$)}, \textcolor{green}{8-mode POD ($\diamond$)}, \textcolor{magenta}{8-mode BPOD ($+$)} (b) First two outputs, symbols as defined in (a).}
  \label{fig:a1b1_KE_growth}
  \end{figure*}

We next compare the impulse response of the system to that of the reduced-order models. \rev{The impulse response of a linear system is important, since the response of the system to any input can be found from the convolution of the impulse response with the input.} Figure~\ref{fig:a1b1_KE_growth} shows the capturing of the growth of kinetic energy by POD and BPOD models, as well as the first two outputs of the reduced-order models. The poor performance of POD at low orders for the traveling structure perturbation is evident. Even the eight-mode POD model, which captures the energy growth well, does not accurately capture the phase of the oscillations.

\begin{figure}%[htpb]
\centering
\includegraphics[width=0.60\linewidth]{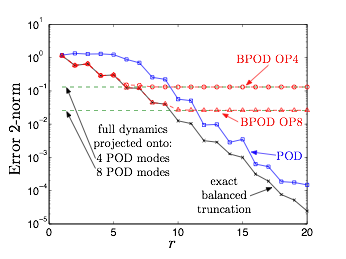}
\caption{Error 2-norms for the $\alpha=1,\beta=1,Re=1000$ perturbation for full balanced truncation, POD and BPOD at two output projections. }
\label{fig:a1b1_twonorms}
\end{figure}

Figure~\ref{fig:a1b1_twonorms} shows the 2-norms of the error between the impulse responses of the reduced-order models and the full simulation $\|G-G_r\|$, normalized by the 2-norm of the impulse response of the full simulation. This figure is a clear demonstration of the effect of output projection. A four-mode output projection means that we are effectively performing balanced truncation on the dynamics of the first four POD modes of the full system. The dashed lines in the figure indicate the 2-norms of the error between the full dynamics of the output-projected system and the full system. As the rank of the BPOD models is increased, the dashed lines, which are the limit of accuracy, are reached fast. As already seen in Fig.~\ref{fig:a1b1_KE_growth}, for very low-order models, standard POD is clearly outperformed by balanced POD. We see that ten POD modes are needed to match the performance of a four-mode output projection BPOD model at the same rank. The slow improvement in POD model performance indicates that the dynamics of the perturbation can not be represented only by retaining the first few most controllable (POD) modes. \sout{Beyond this rank (eight and ten for four- and eight-mode output projections respectively), adding new modes} \rev{Adding new BPOD modes beyond rank eight and ten (for four- and eight-mode output projections, respectively)} does not improve the model performance noticeably, since the dynamics of the output-projected system is already captured fully. It is also important to note that the performance of the BPOD models is identical to that of full balanced truncation almost until the rank at which BPOD model error norms level off due to the output projection. This indicates that the higher balancing modes which are not computed accurately due to the approximation inherent in the output projection (such as those shown in Fig.~\ref{fig:a1b1_balmodes}) do not significantly influence the reduced-order model performance, the main limitation being the capturing of the full system by the output projection.

\subsection{Frequency Response}\label{sec:results-singlefreqresp}

The frequency response encompasses system behavior over the complete range of possible forcing, and is perhaps the best indication of overall system performance. Therefore, from the control designer's point of view, having a low-order model that represents well the frequency response of the original system is of key importance. Frequency response of single-wavenumber perturbations was investigated by Schmid and Henningson~\cite{SchmidHenn} using the resolvent norm, where at each frequency the \sout{maxumum} \rev{maximum} amplification over all initial conditions is computed. Here the frequency response of the system with a given actuator or perturbation is of interest. 

\begin{figure*}%[htpb]
  \centering
\subfigure[]{\includegraphics[width=0.45\linewidth]{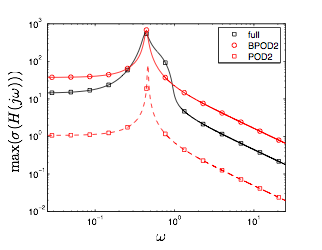}}
\subfigure[]{\includegraphics[width=0.45\linewidth]{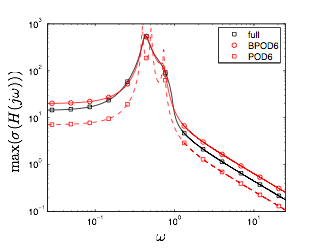}}
\caption {Frequency response of the models for the $\alpha=1,\beta=1$ perturbation. (a) 2-mode POD (dashed), 2-mode BPOD (full), (b) 6-mode POD (dashed), 6-mode BPOD (full).}
\label{fig:a1b1_freqresp}
\end{figure*}

A standard way of representing synthesized frequency response for MIMO (multiple input multiple output) systems is a plot of the maximum singular value of the transfer function matrix $\max(\sigma(H(i\omega)))$ %\remark{Milos, I changed this to $i\omega$ instead of $j\omega$, to be consistent with things like (22)---$j$ is common in electrical engineering, but $i$ is more common in every other field.} 
as a function of frequency, also known as a singular value Bode plot. Fig.~\ref{fig:a1b1_freqresp} (a) and (b) shows such plots for the $\alpha=1,\beta=1$ perturbation and is a clear demonstration of the advantages of BPOD for capturing the dynamics of the system. We see that even for a two-dimensional model the resonant peak is captured well by the model, while for standard POD the peak is very narrow, with very low response at other frequencies. This behavior is typical of balanced truncation, as shown in~\cite{DullerudPag} - the first modes to be captured are the ones which are most significant dynamically, while the correct response is gradually built up in less significant frequency bands as more modes are added. For standard POD models, on the contrary, the response improves incrementally at all frequencies as more modes are added and a higher number of modes is needed to accurately capture the resonant peak. For ten mode models, both standard POD and BPOD capture well the frequency response (not shown in the figure), with BPOD frequency response being almost indistinguishable from the full system one. Standard POD frequency response also includes spurious non-physical peaks at low order, which correspond to eigenvalues very close to the imaginary axis for low order of truncation, as seen in Fig.~\ref{fig:a1b1_spec} (a). 

\begin{figure}
\centering
\includegraphics[width=0.60\linewidth]{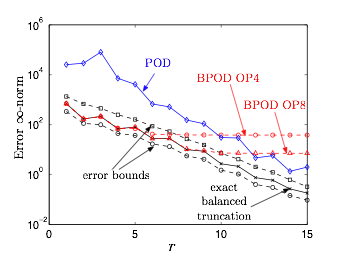}
\caption{Infinity error norms for POD, exact balanced truncation and BPOD with the infinity error bounds.}
\label{fig:a1b1_infnorms}
\end{figure}

We are also interested in the worst-case error between the reduced-order model and the full simulation, which is known as the infinity norm of the system. Balanced truncation has apriori error bounds for the infinity norm~\cite{DullerudPag}. The error $H_\infty$ lower bound for any reduced-order system is
\begin{equation}
 \|G-G_r\|_\infty \ge \sigma_{r+1},
\label{eq:bt_lbound}
\end{equation}
where $G(s)$ is the full model, $G_r(s)$ is the reduced-order model with state dimension~$r$ and $\sigma_j$ is the $j$-th Hankel singular value (in decreasing order). The upper bound for the error is given by
\begin{equation}
 \|G-G_r\|_\infty \le 2\Sigma_{j=r+1}^{n}\sigma_j .
\label{eq:bt_ubound}
\end{equation}
The upper bound on the error can be very close to the lower bound if the HSVs decrease fast. Figure Fig.~\ref{fig:a1b1_infnorms} shows the infinity norm of the error transfer function between the full system and the reduced-order model as a function of model rank for the first fifteen orders of truncation. The infinity norms for exact balanced truncation lie within the theoretical bounds given by Eqs.~(\ref{eq:bt_lbound}) and~(\ref{eq:bt_ubound}), while for BPOD, for each of the two output projections, the norms stay within bounds up to approximately the rank of the output projection, analogous to the two-norms in Fig.~\ref{fig:a1b1_twonorms}. The infinity norms for standard POD at low rank are considerably higher than those for balanced truncation and BPOD, corresponding to the frequency responses shown in Fig.~\ref{fig:a1b1_freqresp} (a) and (b).

\subsection{Variation of Reynolds number}\label{sec:results-singlechangingR}

\begin{figure*}%[htpb]
  \centering
  \includegraphics[width=0.24\linewidth]{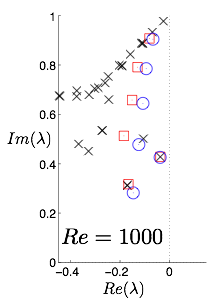}
  \includegraphics[width=0.24\linewidth]{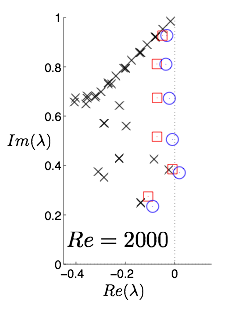}
  \includegraphics[width=0.46\linewidth]{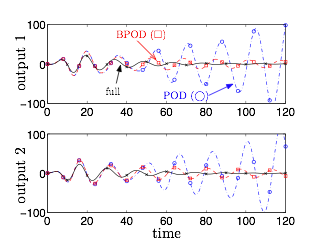}
\caption{Left: The comparison of spectra of the full operator at $\alpha=1,\beta=1$ to the spectra of rank 12 models as the value of $Re$ is increased to 2000. Right: The performance of corresponding reduced-order models at $Re=2000$ - first two outputs. See text for detailed description. Symbols: full ($\times$),  \textcolor{red}{BPOD ($\square$)},  \textcolor{blue}{POD ($\bigcirc$)}.}
\label{fig:a1b1_changeRe}
\end{figure*}

Another very desirable feature of a reduced-order model is good performance for off-design values of the system parameters. We would like the models to remain valid for a wide range of the model parameters, or at least \sout{the appropriate range} \rev{for the range appropriate} for the physical application of the model. The only parameter we are considering in our models is the Reynolds number, so the response of models was compared to the full simulation when $Re$ is changed. Separating the operators from Eq.~(\ref{eq:linearized}) into convective and diffusive parts, we can re-write the state-space equation as
\begin{equation}
\dot{x}=A_{\text{conv}}x + \frac{1}{Re}A_{\text{diff}}x + Bu.
\label{eq:ss-convdiff}
\end{equation}
We can then separately project the matrices $A_\text{conv}$ and $A_\text{diff}$ as in Eq.~(\ref{eq:balanced_system}) at any Reynolds number onto the POD and BPOD modes obtained at $Re=1000$, the $B$ matrix being just the initial condition at $Re=1000$. Figure~\ref{fig:a1b1_changeRe} shows the performance of 12-mode standard POD and BPOD models when the value of $Re$ in Eq.~(\ref{eq:ss-convdiff}) was changed to $2000$ and the impulse response of the resulting models was compared to the impulse response of the full system. This rank of the model was chosen since both models perform well at the design condition of $Re=1000$. We see that the BPOD model eigenvalues stay closer to the full simulation eigenvalues (which move as well), and also remain in the left half of the complex plane, while for $Re=2000$ the standard POD model becomes unstable. This indicates a greater range of validity for BPOD models and better stability at off-design conditions than standard POD. 

\section{Results: three-dimensional localized actuator}
\label{sec:results-3D}

\remark{The name of this section was changed from `3-D disturbance' to `3-D actuator'}

\begin{figure*}
\centering
\includegraphics[width=0.32\linewidth]{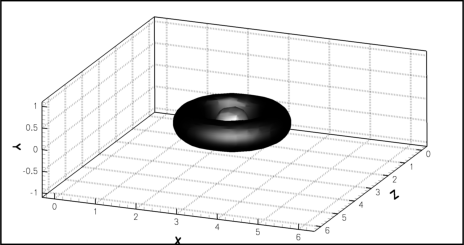}
\includegraphics[width=0.32\linewidth]{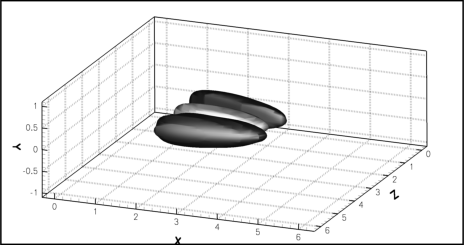}
\includegraphics[width=0.32\linewidth]{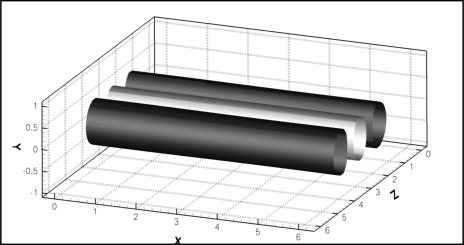}
\caption{The development of the wall-normal velocity of the perturbation given by Eq.~(\ref{eq:blob_IC}) at $t=0$ (left), $t=14$ (middle) and $t=160$ (right) which corresponds to the maximum energy growth. Positive velocity is light, and negative velocity is dark.}
\label{fig:blob_evolution}
\end{figure*}

We next consider a localized body force actuator in the center of the channel which cannot be described by a one-dimensional problem. \rev{This case corresponds to Eq.~(\ref{eq:ss_act_dist}) without the $Fu_2$ term, with the input matrix $B$ representing the velocity field in Fig. \ref{fig:blob_evolution}.} Individual localized disturbances to channel flow were investigated by Henningson et al~\cite{HennLundJoh-93}. Since balanced truncation involves the approximation of the system's Gramians (although in BPOD we do not actually compute the Gramians themselves), we are interested in following both the forward and adjoint impulse-state responses until all transients have completely decayed. The computational box necessary for following individual localized disturbances long enough in time would be very large, and we instead consider a periodic array of small disturbances in the channel center. It should be noted that the behavior of this periodic array can be quite different from the behavior of a single localized disturbance, in particular considering the energy growth, since the periodic array quickly develops into a streamwise-constant vortex. The exact form of the initial condition considered here is \remark{Error in the equation corrected, it was wrong - $(1-r^2)/\alpha^2$ !}
\begin{equation}
v(x,y,z,0)=A\left(1- \frac{r^2}{\alpha^2}\right)e^{(-r^2/\alpha^2 - y^2/\alpha_y^2)}(\cos(\pi y)+1)
\label{eq:blob_IC}
\end{equation}
where $x_c,0,z_c$ are the coordinates of the center of the computational domain and $r^2=(x-x_c)^2+(z-z_c)^2$. The wall-normal vorticity is zero. This form was picked in order to satisfy the condition that the mean perturbation velocity is zero at each wavenumber. The $(\cos(\pi y)+1)$ term was added to make the field satisfy exactly the boundary conditions $v(\pm1)=v_y(\pm1)=0$. The amplitude~$A$ was set to 1 for this simulation, and the parameters $\alpha$ and $\alpha_y$ were set to $\alpha=0.7$ and $\alpha_y=0.6$. The Reynolds number chosen for this simulation was $Re=2000$. The traveling structure rapidly develops into a streamwise-constant structure, since the growth of wall-normal vorticity results in the development of streamwise streaks (see Fig.~\ref{fig:blob_evolution}). 

\begin{figure*}%[htpb]
\centering
\subfigure[]{\includegraphics[width=0.45\linewidth]{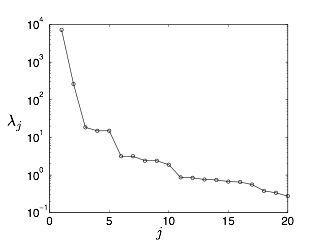}}
\subfigure[]{\includegraphics[width=0.45\linewidth]{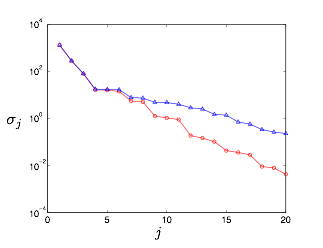}}
\caption{(a) The first 20 POD eigenvalues for the Gaussian-like disturbance impulse response. (b) The first 20 HSVs for five-mode ($\circ$) and ten-mode ($\triangle$) output projections.}
 \label{fig:3Dsigvals}
 \end{figure*}

\begin{figure*}%[htpb]
\centering
\includegraphics[width=0.32\linewidth]{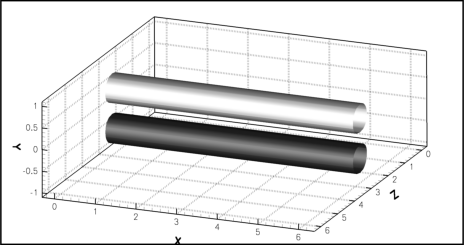}
\includegraphics[width=0.32\linewidth]{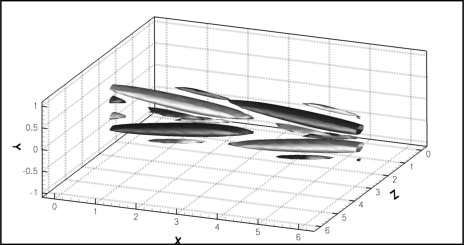}
\includegraphics[width=0.32\linewidth]{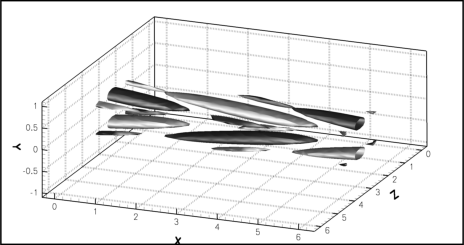}
\caption{The first, fourth and fifth POD modes for the localized actuator, showing streamwise velocity. The iso-surfaces show half of the maximum (light) and minimum (dark) value.}
\label{fig:3DregPOD}
\end{figure*}

The grid size was $32\times 65\times 32$, corresponding to 133,120 states for the full $(v,\eta)$ system. The simulation was ran for 1200 dimensionless time units ($t=t^d U_{c}/\delta$), and the timestep used was $\Delta t=0.004$.  During this time, the energy of the initial disturbance decayed to approximately 1.5 percent of its initial value. The POD modes were taken over 1000 snapshots, with fine spacing between snapshots for the initial period in order to capture the traveling structures well and larger spacing once the streamwise structure was developed, after it was verified that POD eigenvalues and the corresponding modes do not change significantly if more snapshots are used. Fig.~\ref{fig:3Dsigvals} (a) shows the POD eigenvalues of the impulse response. The first five modes contain $99.72\%$ of the perturbation energy, and the first ten modes contain $99.9\%$ of the energy. In this case the spectrum contains both streamwise-constant (and nearly-streamwise constant modes) as well as traveling structures \rev{due to the initial transient}. The first three modes are streamwise-constant structures, while the fourth and the fifth modes correspond to a traveling structure, which accounts for only $0.40\%$ of the total energy. Modes one, four and five are shown in Fig.~\ref{fig:3DregPOD}. 

\begin{figure*}%[htpb]
\centering
\includegraphics[width=0.32\linewidth]{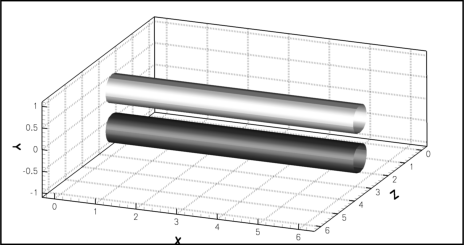}
\includegraphics[width=0.32\linewidth]{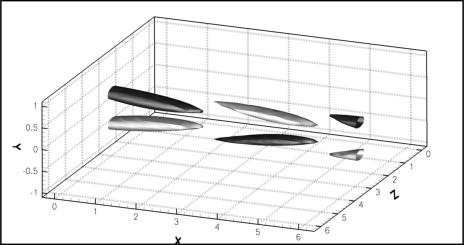}
\includegraphics[width=0.32\linewidth]{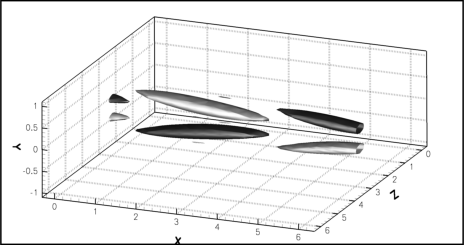}
 \includegraphics[width=0.32\linewidth]{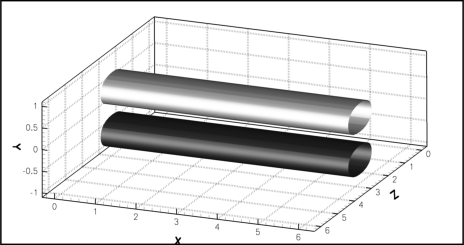}
 \includegraphics[width=0.32\linewidth]{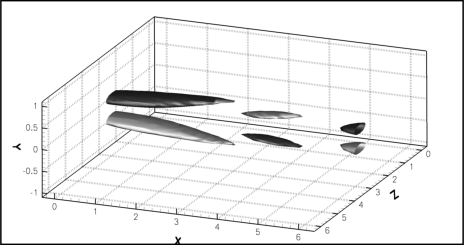}
 \includegraphics[width=0.32\linewidth]{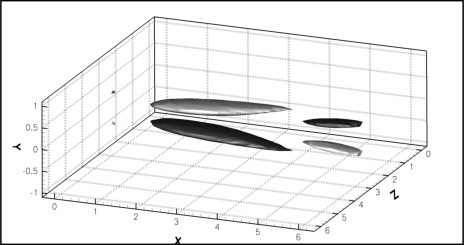}
\caption{Top row: primal modes one, four and five from balanced POD,  showing streamwise velocity for the localized perturbation. The modes are from a five-mode output projection. The iso-surfaces show half of the maximum (light) and minimum (dark) value. Bottom row: the corresponding adjoint modes.  Note the similarity between the primal modes and the corresponding POD modes in Fig.~\ref{fig:3DregPOD}.}
\label{fig:3DbalPOD}
\end{figure*}

Next, the adjoint simulations were ran and the BPOD procedure was performed on a five-mode output projection, containing only the most important traveling structure, as well as on a ten-mode output projection. These ranks were chosen due to large drops in energy significance after the fifth and tenth mode, as shown in Fig.~\ref{fig:3Dsigvals} (a). Fig.~\ref{fig:3Dsigvals} (b) shows the HSVs for these two output projections. We notice that the HSVs are equal for the pairs of modes 4--5 and 7--8 for five-mode output projection, corresponding to traveling structures in the basis of BPOD modes. Even more interestingly, for the ten-mode output projection, HSVs for the modes 4--6 are equal. Although the stability of balanced truncation models is guaranteed only when $\sigma_{r+1}<\sigma_r$, where $r$ is the number of states retained~\cite{DullerudPag}, 4-mode, 5-mode and 7-mode models for both output projections were found to be stable. The model error for impulse response, however, decreases significantly if both modes corresponding to a traveling structure are included, as will be shown in Sec.~\ref{sec:results-3Dimp} (see Fig.~\ref{fig:3Derror2norms}).  Balanced POD modes one, four, and five, and the corresponding adjoint modes for the five-mode output projection are shown in Fig.~\ref{fig:3DbalPOD}. Note that the structure of modes four and five in Fig.~\ref{fig:3DbalPOD} is almost identical, except for a spatial phase shift of exactly one quarter of the periodic domain.

\subsection{Impulse Response}\label{sec:results-3Dimp}

\begin{figure*}%[htbp]
\centering 
\subfigure[]{\includegraphics[width=0.45\linewidth]{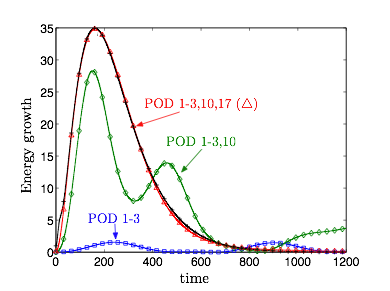}}
\subfigure[]{\includegraphics[width=0.45\linewidth]{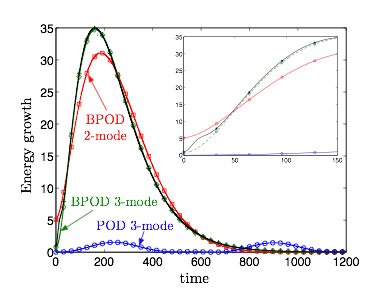}}
\caption{ (a) Three, four and five-mode POD models formed from the indicated modes (b) Two-mode and three-mode BPOD models. The very low-order BPOD models do not capture very well the initial transient, as shown in the inset. \rev{The BPOD models are from the ten-mode output projection.}}
\label{fig:blobimp_resp}
 \end{figure*}

Figure~\ref{fig:blobimp_resp} (a) shows the perturbation energy growth as captured by three different standard POD models. It was observed that the inclusion of modes which come in pairs (see Fig.~\ref{fig:3Dsigvals} (a)) in the basis used to form the reduced order models does not change the system behavior appreciably - the response of a model including modes 1--9 (not shown in figure) is virtually indistinguishable from the response of the model including only the first three modes. Hence, the traveling structure modes do not contribute significantly to the dynamics of this perturbation. The inclusion of the tenth mode, which is streamwise-constant, improves the performance significantly, and the model composed of only the first three modes and the tenth mode performs as well as one including the first ten modes. In the same fashion, including the mode pairs 11--12, 13--14 and 15--16 does not affect the model performance. Including the seventeenth mode, which is also a streamwise-constant mode, improves the performance further. The tenth and the seventeenth mode correspond to $0.025\%$ and $0.0074\%$ of the total energy. The low-order standard POD models were found to capture poorly the initial condition of the full simulation (this will be discussed in more detail in Sec.~\ref{sec:results-singleactuation}), so they were also started from different initial conditions at later times (before or around the peak energy growth), when the projection of the simulation onto POD modes is close to the full simulation data, and they still did not capture the correct peak and the subsequent decay of the energy.

On the other hand, the performance of very low-order BPOD models is significantly better. Figure~\ref{fig:blobimp_resp} (b) shows the perturbation energy growth as captured by three different models. Although the two-mode model does not accurately capture the initial condition, it does represent the energy growth at later stages reasonably well. A three-mode BPOD model captures the kinetic energy of the full simulation very well except for the initial period. While more modes are needed to capture exactly the initial transient, if only the energy growth is of interest, the three-mode model is sufficient. This striking difference is an illustration of the advantage of balanced truncation - for standard POD it is difficult to know apriori which modes will be important for the system dynamics as demonstrated above and a good low-order model was found only after a careful examination of the mode basis which provided some insight into the underlying physics.

\begin{figure}[htpb]
\centering
\includegraphics[width=0.60\linewidth]{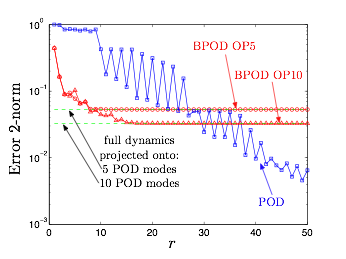}
\caption{Error 2-norms for localized actuator, showing POD models and BPOD at two output projections. } 
\label{fig:3Derror2norms}
\end{figure}

Figure~\ref{fig:3Derror2norms} shows the error $\|G-G_r\|$ for BPOD using the two output projections and for standard POD.  The error for the BPOD models quickly reaches the asymptotic values dictated by the output projection, although more modes are needed compared to the optimal perturbation case described in the previous section due to the more complex dynamics. Standard POD starts to match the performance of the ten-mode output projection BPOD only around rank 30, and also varies a lot with the model rank. This corresponds to the already observed fact that dynamically important POD modes are not highly ranked in terms of energy. Whenever the POD modes come in pairs, including only one of the modes results in deterioration of model performance.  Eventually standard POD has better performance than BPOD, however recall that these POD models have $r$ outputs while the BPOD models have only $s$ outputs. It should also be noted that some POD models exhibit sustained or very slowly decaying oscillations, and that the corresponding two-norms are in fact infinite. Since the simulation time for Fig.~\ref{fig:3Derror2norms} is finite, the two-norms of such models appear to be large but finite as well. Although the four-mode BPOD models are stable for both output projections, including just one of the modes corresponding to a pair of equal HSVs 4--5 for the five-mode output projection deteriorates model performance, while we see a large decrease in the error when the fifth mode is included, as well as when we include subsequent pairs. For the ten-mode output projection, there are three equal HSVs 4--6, and a significant decrease in the error is seen only when we include all three of those modes (in particular, the error norm of the five-mode BPOD model is significantly larger than that of the four-mode model). 

\subsection{Frequency Response}\label{sec:results-3Dfreqresp}

\begin{figure*}%[htpb]
\centering
\subfigure[]{\includegraphics[width=0.45\linewidth]{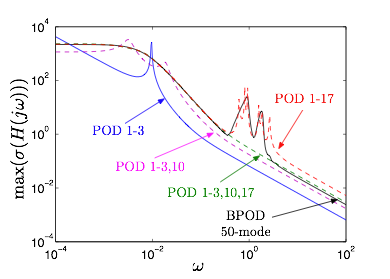}}
\subfigure[]{\includegraphics[width=0.45\linewidth]{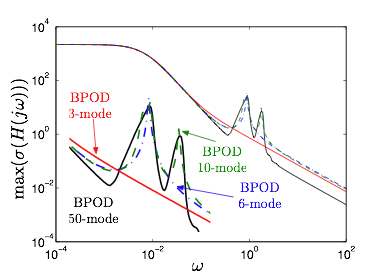}}
\caption{Singular value Bode plots for POD (a) and BPOD (b) models for the localized disturbance. The models are compared to a 50-mode BPOD model. The close-up in (b) shows that a six-mode BPOD model is needed to capture the larger resonant peak, and a ten-mode BPOD model captures both peaks.}
\label{fig:3Dsigmas}
\end{figure*}

Figure \ref{fig:3Dsigmas} shows the singular value Bode plots of standard POD and ten-mode output projection BPOD models for the localized disturbance. The frequency response of the 50-mode BPOD model, which is a very close approximation of the frequency response of the actual disturbance, has the shape of a low-pass filter with a break frequency of $0.01$ $rad/s$ with two resonant peaks near $1$ $rad/s$, which are similar to the peak observed for the single-wavenumber traveling structure perturbation in the previous section. We see that standard POD models again have spurious peaks at low model ranks. The addition of mode pairs corresponding to the traveling structures is necessary in order to reproduce the peaks around $1$ $rad/s$ for both standard POD and BPOD, however BPOD captures those peaks with only the triple 4--6 and the mode pair 7--8, as well as modes 9 and 10 (Figure \ref{fig:3Dsigmas} (b)) while all standard POD modes 1--17 are needed to reproduce the same peaks and there are still spurious peaks. Since the peaks correspond to the low-energy traveling structures, it is not surprising that only a three-mode BPOD model performs so well in capturing the kinetic energy of the full simulation, as shown in the previous section. On the other hand, if the frequency response of the actuator around the frequency of $1$ $rad/s$ needs to be captured accurately, the higher BPOD modes need to be included. 

\subsection{Variation of Reynolds number}\label{sec:results-3DchangingR}

\begin{figure}%[htpb]
\centering
\subfigure[]{\includegraphics[width=0.42\linewidth]{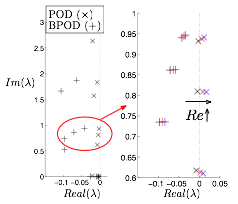}}
\subfigure[]{\includegraphics[width=0.50\linewidth]{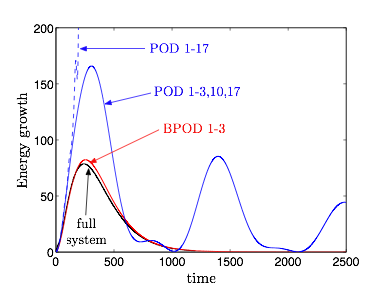}}
\caption{(a) The eigenvalues of 17-mode models at $Re=2000$, $Re=2500$, and $Re=3000$. (b) The performance of the models at the off-design condition of $Re=3000$. \remark{Part (a) was changed to avoid confusion between real part and Reynolds number ($Real(\lambda)$ instead of $Re(\lambda)$)}}
\label{fig:blob_changeRespec}
\end{figure}

Figure~\ref{fig:blob_changeRespec} shows some of the eigenvalues of the 17-mode POD and BPOD models as the Reynolds number is increased. As in Section~\ref{sec:results-singlechangingR}, we use the modes from the design condition ($Re=2000$ in this case) and form the models using Eq.~(\ref{eq:ss-convdiff}). Both POD and BPOD models have eigenvalues on the real axis very close to the origin, which remain stable and correspond to the slow evolution of the streamwise-constant structures.  At each $Re$, the eigenvalues of both models move towards the right half of the complex plane and while the BPOD model always remains stable, the standard POD model first appears marginally stable at $Re=2500$ and then unstable at $Re=3000$. The effect of the eigenvalues that move to the right half of the complex plane is clearly seen in Fig.~\ref{fig:blob_changeRespec} (b). A model that includes modes 1--17 grows unstable quickly at $Re=3000$, showing that inclusion of modes which at design condition do not contribute signifcantly to the overall dynamics can significantly deteriorate the performance of the model at off-design condition. This can also be seen from the frequency responses shown in the previous section - even at design condition, the spurious high peaks correspond to marginally stable modes. Although stable, the 1,2,3,10,17 POD model is highly inaccurate at $Re=3000$, with large peaks in the kinetic energy which decay very slowly, indicating the high sensitivity of those POD models which remain stable to a change in the Reynolds number. On the other hand, the three-mode BPOD model is still remarkably close to the full system. 

It is important to note here that as the Reynolds number is increased, the \sout{non-linearity} \rev{nonlinearity} will have a stronger effect on the development of the disturbance and the reduced-order model may not be valid for higher $Re$ in the first place. The comparison of the linear perturbation growth with a full nonlinear DNS solution is essential for a true validation of the models for controls applications, since we may be modeling the linearized flow well, but the linearized flow may not be a good approximation to the actual flow. This comparison is subject of current work. 

\subsection{Capturing of actuation}\label{sec:results-singleactuation}

\begin{figure*}[htpb]
\centering
\subfigure[]{\includegraphics[width=0.46\linewidth]{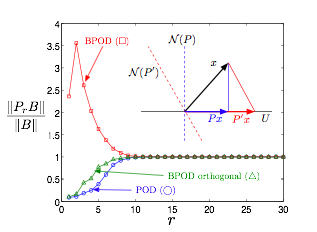}}
\subfigure[]{\includegraphics[width=0.46\linewidth]{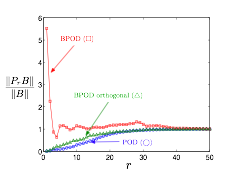}}
\caption{Norm of the projection of the $B$ matrix (a single column vector) onto subspaces used for reduced-order models, (a) $\alpha=1,\beta=1$, (b) localized disturbance. The diagram in (a) illustrates the non-orthogonal projection used in BPOD. \remark{Error in plot fixed - the green curves now say 'BPOD orthogonal', not 'POD orthogonal'}}
\label{fig:proj_on_b}
\end{figure*}

An important property of a reduced-order model is how well it captures the effects of the actuator in the original system, especially for models that are intended for developing controllers. In order for a reduced-order model to capture the effect of an actuator, it is necessary at a minimum for the input term in the equations ($Bu$ in Eq.~(\ref{eq:primal})) to be contained in the subspace used for projecting the equations. Note that here, even for the standard POD case, the effect of the ``actuator'' is partially included, since the dataset used for POD is generated by an impulsive input. One way to measure the degree to which the input ``directions'' are captured by the modes used in the model is to compute the projection of the columns of the input matrix $B$ in Eq.~(\ref{eq:primal}) onto the basis modes. In the system we are considering here, $B$ is a single column vector, representing the initial disturbance given to the system (or actuation via a body force in the center of the domain). Fig.~\ref{fig:proj_on_b} shows the norm of the projection $\|P_rB\|/\| B\|$ of the standard and BPOD modes onto the input vector $B$, which is just the initial condition for each simulation. The balancing modes clearly capture the input direction with many fewer modes than standard POD: even very low-order models have a significant norm after projection, and in fact the norm of~$B$ after projection is almost always greater than the norm of $B$ due to the non-orthogonal projection, as shown in Fig.~\ref{fig:proj_on_b}. Any orthogonal projection such as $P$ must satisfy $\|Px\|\le\|x\|$, while for a non-orthogonal projection we may have $\|P' x\|>\|x\|$, which is the case for the first several BPOD modes. Clearly, the $B$ matrix has a very small projection onto the POD modes for standard POD unless many modes are taken, so it is impossible for very low-order POD models to capture the response of an actuator without introducing more modes (such as the $B$ matrix itself, Krylov subspaces, or shift modes~\cite{Noack-03}). 

\subsection{Subspace comparison}\label{sec:subspace_comp}

\begin{figure}%[htpb]
\centering
\includegraphics[width=0.60\linewidth]{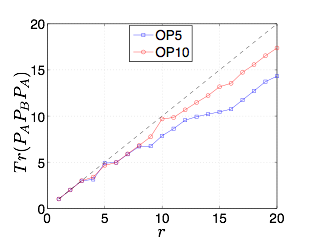}
\caption{Plots of $Tr(P_AP_BP_A)$ as a function of the subspace rank $r$ for the localized perturbation for both output projections.}
\label{fig:traces_blob}
\end{figure}

The BPOD procedure uses both a different projection and a different set of modes in order to form reduced-order models, and we next look at a comparison of the two mode subspaces. A way to compare two subspaces is to compute $Tr(P_AP_BP_A)=T$, where $P_A$ and $P_B$ are the corresponding projection operators~\cite{EveSir-95}. The trace of the matrix $Tr(P_AP_BP_A)$ as a function of the subspace rank is shown in Fig.~\ref{fig:traces_blob} for the five-mode and ten-mode output projections where $P_A$ and $P_B$ are the orthogonal projectors onto the POD and BPOD subspaces respectively. The value of the trace $T$ %\todo{Milos, why $s$ here when you just called this $T$?}  
is the same as the subspace rank $r$ at low order, indicating very similar modes (see Figs.~\ref{fig:3DregPOD} and \ref{fig:3DbalPOD}). For both POD and five-mode OP BPOD, modes four and five are a pair of structures, and including both modes from the pair brings the value of $s$ to almost exactly 5. It is interesting to observe that for the five-mode output projection, $r=T$ exactly at $r=5$, while for the ten-mode output projection $r=T$ at $r=10$, and that above those values the value of the trace is lower than the rank. This can be explained by the fact that BPOD is attempting to approximate the output projection of the data of the given rank. It is interesting to note that the subspaces including the first three POD and BPOD modes are virtually identical, indicating that the non-orthogonal Petrov-Galerkin projection via adjoint modes makes the enormous difference that we have seen in the performance of the corresponding models. As mentioned above, the POD basis is the basis of the most controllable modes, and is indeed optimal in capturing a given dataset, but as we have shown, it can fail to capture the dynamics correctly. 

\section{Conclusions and Further Work}\label{sec:conclusions}
\remark{We have added a few sentences here to address the concerns of the reviewers.}

We have shown that the balanced POD procedure offers several advantages over the standard POD technique for modeling linearized channel flow. Balanced POD makes the benefits of the systematic model reduction via balanced truncation available for large systems, resulting in models that capture the impulse response, frequency response and the effects of an actuator much better than their corresponding POD models.  The differences between POD and BPOD are especially striking for localized actuators or disturbances which are often of interest in applications. As BPOD allows model reduction for very large systems, the localized perturbation models were obtained using three-dimensional fields in physical space, without the standard wavenumber decoupling which can be unwieldy for localized disturbances, \rev{since model reduction would have to be performed at each wavenumber pair resulting in models which are still of relatively high rank, and each wavenumber pair would need to be controlled separately. More importantly, this approach allows the extraction of dynamically dominant structures in physical space which may have contributions from many wavenumbers. This feature of the method will be very important in further applications of this method, and in particular in spatially developing flows.}  

There is an indication that the non-orthogonal projection onto the balancing modes with the use of adjoint modes plays a key role in the model performance.  The subspaces spanned by the BPOD modes are very close to those spanned by the POD modes (which are optimal in representing the impulse response dataset).  However, the dynamical models produced by BPOD are quite different, because of the non-orthogonal projection (or, equivalently, an orthogonal projection with a new inner product, defined by the observability Gramian~\cite{Rowley-ijbc05}).

This work has focused on the linearized system in order to characterize the BPOD model reduction method. A true test of any controller is the application to the full nonlinear simulation, and verification of reduced perturbation growth, and this is a subject of our ongoing work. \rev{In addition, several techniques are available for obtaining nonlinear reduced-order models.  One method for snapshot-based balancing for nonlinear systems has been introduced by Lall et al~\cite{LallMarsden-02}, but this involves considerably more computation than the present method, and is not feasible for large systems.  A simpler, but presumably less effective approach is to project the full Navier-Stokes equations onto the BPOD modes computed for a linearized system; such a procedure involves no additional computational expense over the methods presented in this paper.}

Although the body force actuation we discuss here would be desirable in applications, it is not as practically feasible as wall blowing and suction, which has been studied extensively and has shown promise both in computations and experiment~\cite{HogbergBewley-03,Lundell-03}. \rev{We also note that blowing and suction could be incorporated into the control framework described here in a relatively straightforward fashion, via the lifting approach described in Ref.~\cite{HogbergBewley-03}.} \sout{Another topic of our current work is the inclusion of wall blowing/suction actuators in the linearized system.} Finally, the BPOD procedure is applicable to other flows that can be linearized, such as Couette flow, boundary layer flows, or flows over airfoils. Moreover, the method can be applied to any large linear system and can therefore be a useful addition to the tools of modern control theory.

\subsection*{Acknowledgments}
We gratefully acknowledge Dan Henningson, Shervin Bagheri, Espen \AA kervik and Sunil Ahuja, as well as the anonymous referees, for helpful comments and suggestions.  This work was supported by the National Science Foundation, grant CMS-0347239.

\appendix
%\section{Independence of Balanced Truncation from Inner Product}\label{sec:App-indep_ip}
\section{Snapshot-based balanced truncation using a continuous adjoint}\label{sec:App-indep_ip}

\remark{The title and the contents of this appendix were changed considerably in order to clarify the details of computing BPOD in the case when a continuous adjoint is used}
%\todor{Some things were not entirely clear in the appendix, making some changes based on a discussion with Sunil.}

When computing the exact balanced truncation, the balancing transformation is found from the eigenvalue problem $W_cW_oT = T\Sigma^2$ where $W_c$ and $W_o$ are the controllability and observability Gramians with $A^+=A^T$. We show here that, although the product $W_cW_o$ does not depend on the inner product on the state space used to define the adjoint system, the appropriate weight $M$ needs to be included in the computation via the method of snapshots.
 
We can represent the weighted inner product of two vectors $q_1$ and $q_2$ as 
 %
 %\begin{equation}
 %\ip<q_1,q_2>_M=q_2^TMq_1
 %\end{equation}
  \begin{equation}
\ip<x_1,x_2>_M = \int_\Omega x_1^*Mx_2d\Omega
\end{equation}
 where the domain of integration $\Omega$ is the Hilbert space itself. The star denotes the complex conjugate \rev{transpose}. The inner product weight~$M$ is part of the definition of the Hilbert space itself. We define the so-called continuous adjoint of an operator $A$ with respect to this inner product as
 \begin{equation}\label{eqn:mipdef}
 \ip<Aq_1,q_2>_M=\ip<q_1,A^+q_2>_M
 \end{equation}
We use the symbol $+$ in order to distinguish the adjoint from the standard matrix transpose $A^T$. From this definition it is easily shown that 
 \begin{eqnarray}\label{adjoints}
 \begin{aligned}
 \ip<Ax,z>_M = \ip<x,A^+z>_M \quad & \Rightarrow \quad A^+ = M^{-1}A^TM   \\
 \ip<Bu,x>_M = \ip<u,B^+x>  \quad & \Rightarrow \quad B^+=B^TM  \\
 \ip<Cx,y> = \ip<x,C^+y>_M  \quad & \Rightarrow \quad C^+ =M^{-1}C^T  \\
 \end{aligned}
 \end{eqnarray}
In the above, we have assumed that the input and output spaces use the standard (unweighted) inner product.  Next, we obtain for the Gramians:
 \begin{eqnarray}
 \begin{aligned}
 G_c &=\int_{0}^{\infty}e^{At}BB^+e^{A^+t}dt=\int_{0}^{\infty}e^{At}BB^TMM^{-1}e^{A^Tt}Mdt\\
 G_o &=\int_{0}^{\infty}e^{A^+t}C^+Ce^{At}dt=\int_{0}^{\infty}M^{-1}e^{A^Tt}MM^{-1}C^TCe^{At}dt
 \end{aligned}
 \label{eq:modgrams}
 \end{eqnarray}
 where $G_c$ and $G_o$ denote the Gramians obtained with the weighted inner product. Since the matrices $M$ and $M^{-1}$ are constant, we can take them out of the integrals, obtaining
 \begin{equation}
 G_cG_o=W_cW_o
 \end{equation}
 
Thus, we have shown that balanced truncation does not depend on the choice of the inner product used to derive the adjoint system, and this allows us to use a convenient inner product.  (In numerical simulations the `simple' discrete adjoint $A^+=A^T$ may in fact be more difficult to compute than a continuous adjoint which may retain a similar form of the equations; for instance, this is the case for linearized channel flow).
 
%On the other hand, weight \sout{matrix} \todor{matrices for the input and output space inner products do} \sout{does} not cancel \sout{when using different inner products to derive $B^+$ or $C^+$} \todor{in the product $G_cG_o$} (since each appears only in one of the Gramians), so that in this case the choice of inner product does affect \sout{$G_cG_o$} \todor{the balancing transformation}.

Next, we consider the computation of balancing and adjoint modes via the method of snapshots. From the definition of the empirical Gramians (Eq.~(\ref{eq:empgrams})) it is easily shown that $Y^+=Y^TM$ (recall that the snapshots of the adjoint simulations, which are the columns of $Y$, are given by $z(t)=e^{A^+t}C^+$). Thus, we can write the SVD in Eq.~(\ref{eq:SVD}) as
\begin{equation}
Y^TMX=U\Sigma V^T
\end{equation}
%
%One way to form the balancing transformation would be to take $\Phi=XV\Sigma^{-1/2}$ and $\Psi^T =  \Sigma^{-1/2}U^TY_M^TM= \Sigma^{-1/2}U^TY_D^T$. In that case, both the balancing modes and the adjoint modes would be the same as in the case of the continuous adjoint. However, the computation of $\Psi$ may not be possible when using the method of snapshots because the matrix $M$ may be too large or we only have $M$ as a continuous operator. 
If we define the inverse of the balancing transformation as $\Psi^T=\Sigma^{-1/2}U^TY^T$ we can easily compute the adjoint modes just from the SVD and from the adjoint snapshots. Recall that the  \textit{columns} of $\Psi$ give the adjoint modes. The two sets of modes will now be bi-orthogonal with respect to the $M$ inner product, so that $\Psi^TM\Phi=I$. 

An alternative, more intuitive explanation is that, since both the direct and the adjoint snapshots `live' in the state space, the correct inner product is that including the weight $M$ (which is a part of the definition of the Hilbert space in which they reside).  It is therefore this weighted inner product that should be used for forming the matrix for the SVD.  Furthermore, the balancing and adjoint modes are bi-orthogonal with respect to this weighted inner product.

\pagebreak

\bibliography{jabbrv,master}

\begin{thebibliography}{10}

\bibitem{Aakervik_et_al-jfm-07}
E.~\AA{}kervik, J.~H\oe{}pffner, U.~Ehrenstein, and D.~S. Henningson.
\newblock Optimal growth, model reduction and control in a separated
  boundary-layer flow using global eigenmodes.
\newblock {\em J. Fluid Mech.}, 579:305--314, 2007.

\bibitem{Bamieh-01}
B.~Bamieh and M.~Dahleh.
\newblock Energy amplification in channel flows with stochastic excitation.
\newblock {\em Phys.\ Fluids}, 13(11):3258--3269, Nov. 2001.

\bibitem{ButlerFarrell-92}
K.~Butler and B.~Farrell.
\newblock Three-dimensional optimal perturbations in viscous shear flow.
\newblock {\em Phys. Fluids A}, 4:1637--1650, 1992.

\bibitem{DullerudPag}
G.~E. Dullerud and F.~Paganini.
\newblock {\em A Course in Robust Control Theory: A Convex Approach}, volume~36
  of {\em Texts in Applied Mathematics}.
\newblock Springer-Verlag, 1999.

\bibitem{EveSir-95}
R.~Everson and L.~Sirovich.
\newblock Karhunen-{Lo\`eve} procedure for gappy data.
\newblock {\em J. Opt. Soc. Am. A}, 12(8):1657--1664, Aug. 1995.

\bibitem{Farrell-88}
B.~Farrell.
\newblock Optimal excitation of perturbations in viscous shear flow.
\newblock {\em Phys. Fluids}, 31:2093--2102, 1988.

\bibitem{Farrell-01}
B.~F. Farrell and P.~J. Ioannou.
\newblock Accurate low-dimensional approximation of the linear dynamics of
  fluid flow.
\newblock {\em J. Atmospheric Sci.}, 58:2771--2789, 2001.

\bibitem{Gustavsson-86}
L.~H. Gustavsson.
\newblock Excitation of direct resonances in plane poiseuille flow.
\newblock {\em Stud. Appl. Math.}, 75:227--248, 1986.

\bibitem{Hanifi-96}
A.~Hanifi, P.~J. Schmid, and D.~S. Henningson.
\newblock Transient growth in compressible boundary layer flow.
\newblock {\em Phys. Fluids}, 8(3), March 1996.

\bibitem{HennLundJoh-93}
D.~S. Henningson, A.~Lundbladh, and A.~V. Johansson.
\newblock A mechanism for bypass transition from localized disturbances in
  wall-bounded shear flows.
\newblock {\em J. Fluid Mech.}, 250:169--207, 1993.

\bibitem{HennSchmid-92}
D.~S. Henningson and P.~J. Schmid.
\newblock Vector eigenfunction expansions for plane channel flows.
\newblock {\em Stud. Appl. Math.}, 87:15--43, 1992.

\bibitem{HogbergBewley-03}
M.~{H\"ogberg}, T.~R. Bewley, and D.~S. Henningson.
\newblock Linear feedback control and estimation of transition in plane channel
  flow.
\newblock {\em J. Fluid Mech.}, 481:149--175, 2003.

\bibitem{HLB-96}
P.~Holmes, J.~L. Lumley, and G.~Berkooz.
\newblock {\em Turbulence, Coherent Structures, Dynamical Systems and
  Symmetry}.
\newblock Cambridge University Press, Cambridge, UK, 1996.

\bibitem{IlaRow-06}
M.~Ilak and C.~W. Rowley.
\newblock Reduced-order modeling of channel flow using traveling {POD} and
  balanced {POD}.
\newblock AIAA Paper 2006-3194, 3rd AIAA Flow Control Conference, June 2006.

\bibitem{JosSpeKim-97}
S.~S. Joshi, J.~L. Speyer, and J.~Kim.
\newblock A systems theory approach to the feedback stabilization of
  infinitesimal and finite-amplitude disturbances in plane poiseuille flow.
\newblock {\em J. Fluid Mech.}, 332:157--184, 1997.

\bibitem{JovBam-05}
M.~R. Jovanovi\'c and B.~Bamieh.
\newblock Componentwise energy amplification in channel flows.
\newblock {\em J. Fluid Mech.}, 534:145--183, 2005.

\bibitem{KimMoinMoser-87}
J.~Kim, P.~Moin, and R.~Moser.
\newblock Turbulence statistics in fully-developed channel flow at low
  {Reynolds} number.
\newblock {\em J. Fluid Mech.}, 177:133--166, Apr. 1987.

\bibitem{LallMarsden-02}
S.~Lall, J.~E. Marsden, and S.~Glava\v{s}ki.
\newblock A subspace approach to balanced truncation for model reduction of
  nonlinear control systems.
\newblock {\em Int. J. Robust Nonlinear Control}, 12:519--535, 2002.

\bibitem{Laub-87}
A.~J. Laub, M.~T. Heath, C.~C. Page, and R.~C. Ward.
\newblock Computation of system balancing transformations and other
  applications of simultaneous diagonalization algorithms.
\newblock {\em IEEE Trans.\ Automat.\ Contr.}, 32:115--122, 1987.

\bibitem{LeeCortelezzi-01}
K.~H. Lee, L.~Cortelezzi, J.~Kim, and J.~Speyer.
\newblock Application of reduced-order controller to turbulent flows for drag
  reduction.
\newblock {\em Phys.\ Fluids}, 13(5):1321--1330, 2001.

\bibitem{Lundell-03}
F.~Lundell.
\newblock Pulse-width modulated blowing/suction as a flow control actuator.
\newblock {\em Experiments in Fluids}, 35:502--504, 2003.

\bibitem{Min-06}
T.~Min, S.~M. Kang, J.~L. Speyer, and J.~Kim.
\newblock Sustained sub-laminar drag in a fully developed channel flow.
\newblock {\em J. Fluid Mech.}, 558:309--318, 2006.

\bibitem{MoinMoser-89}
P.~Moin and R.~Moser.
\newblock Characteristic-eddy decomposition of turbulence in a channel.
\newblock {\em J. Fluid Mech.}, 200:471--509, 1989.

\bibitem{Moore-81}
B.~C. Moore.
\newblock Principal component analysis in linear systems: Controllability,
  observability, and model reduction.
\newblock {\em IEEE Trans.\ Automat.\ Contr.}, 26(1):17--32, Feb. 1981.

\bibitem{Noack-03}
B.~Noack, K.~Afanasiev, M.~Morzy\'nski, G.~Tadmor, and F.~Thiele.
\newblock A hierarchy of low-dimensional models for the transient and
  post-transient cylinder wake.
\newblock {\em J. Fluid Mech.}, 497:335--363, 2003.

\bibitem{ObinataAnderson}
G.~Obinata and B.~D.~O. Anderson.
\newblock {\em Model Reduction for Control System Design}.
\newblock Springer-Verlag, 2000.

\bibitem{Orszag-1971}
S.~A. Orszag.
\newblock Accurate solution of the orr-sommerfeld stability equation.
\newblock {\em J. Fluid Mech.}, 50:689--703, 1971.

\bibitem{ReddyHenn-93}
S.~C. Reddy and D.~S. Henningson.
\newblock Energy growth in viscous channel flows.
\newblock {\em J. Fluid Mech.}, 252:209--238, 1993.

\bibitem{Reddy-98}
S.~C. Reddy, P.~J. Schmid, J.~S. Baggett, and D.~S. Henningson.
\newblock On stability of streamwise streaks and transition thresholds in plane
  channel flows.
\newblock {\em J. Fluid Mech.}, 365:269--303, 1998.

\bibitem{Rowley-ijbc05}
C.~W. Rowley.
\newblock Model reduction for fluids using balanced proper orthogonal
  decomposition.
\newblock {\em Int. J. Bifurcation Chaos}, 15(3):997--1013, Mar. 2005.

\bibitem{Rowley-physd04}
C.~W. Rowley, T.~Colonius, and R.~M. Murray.
\newblock Model reduction for compressible flows using {POD} and {Galerkin}
  projection.
\newblock {\em Phys.\ D}, 189(1--2):115--129, Feb. 2004.

\bibitem{RowlMars-00}
C.~W. Rowley and J.~E. Marsden.
\newblock Reconstruction equations and the {Karhunen}-{Lo{\`e}ve} expansion for
  systems with symmetry.
\newblock {\em Phys.\ D}, 142:1--19, Aug. 2000.

\bibitem{SchmidHenn}
P.~J. Schmid and D.~S. Henningson.
\newblock {\em Stability and Transition in Shear Flows}, volume 142 of {\em
  Applied Mathematical Sciences}.
\newblock Springer-Verlag, 2001.

\bibitem{Sirovich-87}
L.~Sirovich.
\newblock Turbulence and the dynamics of coherent structures, parts {I--III}.
\newblock {\em Q.\ Appl.\ Math.}, XLV(3):561--590, Oct. 1987.

\bibitem{Smith-thesis}
T.~R. Smith.
\newblock {\em Low-dimensional models of plane {Couette} flow using the proper
  orthogonal decomposition}.
\newblock PhD thesis, Princeton University, 2003.

\bibitem{Tadmor2004tcst2vort}
G.~Tadmor.
\newblock Observers and feedback control for a rotating vortex pair.
\newblock {\em IEEE Trans.\ Contr.\ Sys.\ Tech.}, 12:36--51, 2004.

\bibitem{Trefethen-93}
L.~N. Trefethen, A.~E. Trefethen, S.~C. Reddy, and T.~A. Driscoll.
\newblock Hydrodynamic stability without eigenvalues.
\newblock {\em Science}, 261:578--584, July 1993.

\bibitem{Willcox-02}
K.~Willcox and J.~Peraire.
\newblock Balanced model reduction via the proper orthogonal decomposition.
\newblock {\em AIAA J.}, 40(11):2323--2330, Nov. 2002.

\bibitem{ZhoDoy-98}
K.~Zhou and J.~C. Doyle.
\newblock {\em Essentials of Robust Control}.
\newblock Prentice-Hall, 1998.

\end{thebibliography}
\end{document}